# Asymptotic normality for the counting process of weak records and δ-records in discrete models

RAÚL GOUET[1], F. JAVIER LÓPEZ[2],[*] and GERARDO SANZ[2],[**]

[1]*Departameto Ingeniería Matemática and Centro de Modelamiento Matemático, Universidad de Chile, Casilla 170/3, Correo 3, Santiago, Chile.* E-mail: *rgouet@dim.uchile.cl*
[2]*Departamento Métodos Estadísticos, Facultad de Ciencias, Universidad de Zaragoza, C/Pedro Cerbuna, 12, 50009 Zaragoza, Spain.* E-mail: [*]*javier.lopez@unizar.es;* [**]*gerardo.sanz@unizar.es*

Let $\{X_n, n \geq 1\}$ be a sequence of independent and identically distributed random variables, taking non-negative integer values, and call $X_n$ a $\delta$-record if $X_n > \max\{X_1, \ldots, X_{n-1}\} + \delta$, where $\delta$ is an integer constant. We use martingale arguments to show that the counting process of $\delta$-records among the first $n$ observations, suitably centered and scaled, is asymptotically normally distributed for $\delta \neq 0$. In particular, taking $\delta = -1$ we obtain a central limit theorem for the number of weak records.

*Keywords:* Central limit theorem; martingale; record; weak record

## 1. Introduction

The theory of records is a well established branch of extreme value theory with interesting results from both a theoretical and a practical point of view. See the books by Ahsanullah [1], Arnold *et al.* [2] or Nevzorov [18] for the theory and applications of record and record-related statistics. Given a sequence $\{X_n, n \geq 1\}$ of random variables, an observation $X_i$ is called a record if it is greater than all previous observations; that is, writing $M_n$ for the maximum of the $n$ first observations, if $X_i > M_{i-1}$. If the random variables $X_n$ are integer-valued, an observation is called a weak record if it is greater than or equal to the previous maximum; that is, if $X_i \geq M_{i-1}$ or, equivalently, $X_i > M_{i-1} - 1$. This leads us to consider the following natural extension of the concept of records: for $\delta \in \mathbb{R}$, an observation $X_i$ is called a $\delta$-record if $X_i > M_{i-1} + \delta$, that is, if it is greater than the previous maximum plus a (negative or positive) fixed value $\delta$. For $\delta < 0$, every record is a $\delta$-record, while for $\delta > 0$ this is not the case. Usual records are obtained by taking $\delta = 0$ and, for integer-valued random variables, $\delta = -1$ yields weak records. In this paper we focus attention on the process $N_n^\delta = \sum_{i=1}^n \mathbf{1}_{\{X_i > M_{i-1}+\delta\}}$, counting the number of $\delta$-records among the first $n$ observations, where $\mathbf{1}_{\{\cdot\}}$ stands for the indicator function. An arbitrary value can be given to $M_0$ because we are dealing with asymptotic results.

In addition to being a natural generalization of records and weak records, our concept of $\delta$-record and the study of the associated counting process $N_n^\delta$ can be relevant, among other things, in insurance applications, where one is interested not only in record claims, but also in claims that are close to being records; see, for instance, Balakrishnan *et al.* [5], Hashorva [12] or Hashorva and Hüsler [13]. In fact, the study of observations near the maximum has attracted much attention in the past years, both in the case of fixed size samples (Li [16]; Pakes [19]; Pakes and Steutel [20]) and when observations are considered sequentially (Balakrishnan *et al.* [4, 5] and Khmaladze *et al.* [15]), where we find concepts closely related to $\delta$-records defined in the present work. Khmaladze *et al.* [15] defined the $\varepsilon$-repeated records as the observations $X_i$ which fall in the interval $(M_n - \varepsilon, M_n]$ for $i$ ranging from $\tau_n = \inf\{k : X_k = M_n\}$ (the moment when the maximum $M_n$ is attained) to $n$. Khmaladze's process $Z_n$, counting $\varepsilon$-repeated records, and our $N_n^\delta$ are related by the equation $Z_n = N_n^\delta - N_{\tau_n-1}^\delta$, with $\delta = -\varepsilon$. In Balakrishnan and Stepanov [6] and Khmaladze *et al.* [15], the asymptotic behaviour of $Z_n$ for sequences of independent identically distributed continuous random variables is studied. On the other hand, Balakrishnan *et al.* [5] defined, for fixed $a > 0$, the near-$n$th records as observations $X_i$ in $(X(n) - a, X(n)]$ for $i \in (L(n), L(n+1))$, where $L(n)$ is the $n$th record time and $X(n)$ is the $n$th record value. The number $\xi_n(a)$ of Balakrishnan's near-$n$th records is related to the number of $\delta$-records through $N_{L(n)}^\delta = \sum_{k=1}^n \xi_k(a) + n$, with $\delta = -a$. The asymptotic behaviour of the number of near-$n$th records is considered in that paper for sequences of independent and identically distributed continuous random variables. Finally, we mention $\delta$-exceedance records, defined in Balakrishnan *et al.* [4] for $\delta > 0$, as observations that exceed the previous $\delta$-exceedance by at least $\delta$; in other words, if $X_{T_k}$ is the $k$th exceedance, the following one is $X_{T_{k+1}}$, with $T_{k+1} = \min\{j > T_k | X_j > X_{T_k} + \delta\}$. Clearly, $\delta$ exceedances and $\delta$-records are not equivalent concepts, because for $\delta > 0$, a $\delta$-record is always a $\delta$ exceedance but not conversely.

The behaviour of the number of usual records $N_n^0$ is well understood when the underlying variables $X_n$ are independent and identically distributed with continuous distribution function because, as shown in Renyi [21], the indicators $I_n = \mathbf{1}_{\{X_n > M_{n-1}\}}$ are independent, with $E(I_n) = 1/n$ and, consequently, many asymptotic results for $N_n^0$ are readily obtained. The study of records and weak records in discrete distributions, where the independence of indicators is lost, was initiated by Vervaat [22]. Asymptotic results for the number of records and weak records, including a central limit theorem, for the geometric distribution have been obtained by Bai *et al.* [3]. Strong laws of large numbers and central limit theorems for $N_n^0$ were given by Gouet *et al.* [9, 10] for large classes of discrete distributions classified in terms of their discrete failure rates. See also Key [14] for a law of large numbers for weak records in heavy-tailed discrete distributions.

In this work we obtain central limit theorems for the number of $\delta$-records $N_n^\delta$, $\delta \neq 0$, when the random variables $X_n$ are independent and identically distributed with discrete distribution function $F$ on the non-negative integers. As a particular case, taking $\delta = -1$, we obtain a central limit theorem for the number of weak records. To the best of the authors' knowledge, all the results in this paper are new for $\delta \neq -1$; for $\delta = -1$, they greatly extend the known results for the geometric distribution to a wide class of discrete models.

Our proofs are based on a martingale approach whereby the counting process $N_n^\delta$ is centered by a non-predictable process built from what we call discrete $\delta$ failure rates [see (2.1)]. Asymptotic normality is established using a martingale central limit theorem, requiring the convergence of conditional variances and a Lyapunov-type condition. Both convergence problems are reduced to the study of partial sums of minima of independent identically distributed random variables, whose asymptotic behaviour has been investigated in detail, especially by Deheuvels [7]. Martingales have already proved to be useful in the study of extremes in discrete settings; see Gouet *et al.* [9, 10].

Here we do not consider the case of continuous distributions, unlike the above cited works on recordlike statistics (Balakrishnan *et al.* [4, 5]; Balakrishnan and Stepanov [6]; Khmaladze *et al.* [15]), which were concerned only with continuous distributions. The study of $N_n^\delta$ in the continuous distribution setting is far from trivial for $\delta \neq 0$, because indicators $\mathbf{1}_{\{X_n > M_{n-1} + \delta\}}$ are neither independent nor distribution-free (see Remark 2.1). We center here on integer valued random variables, thus including the especially interesting case of weak records.

The structure of the paper is as follows. Section 2 presents the notation and three preliminary results. The central limit theorems for the number of $\delta$-records, for $\delta < 0$ and $\delta > 0$, are shown in Sections 3 and 4, respectively. Section 5 is devoted to the application of our results to well-known discrete distributions. Finally, the martingale central limit theorem and Deheuvels' theorem on sums of partial minima are presented in the Appendix.

## 2. Notation and preliminary results

Let $\{X_n, n \geq 1\}$ be a sequence of non-negative, integer-valued, independent and identically distributed random variables, with common distribution function $F$, such that $P[X_n = k] = p_k > 0$ for $k \in \mathbb{Z}_+ = \{0, 1, \ldots\}$ and $n \geq 1$ ($p_m = 0$ for $m \leq -1$). Clearly then, $\inf\{x|F(x) \geq 1\} = \infty$. The inverse of any distribution function, say $G$, will be denoted $G^-(y) = \inf\{x|G(x) \geq y\}$ for $0 \leq y \leq 1$.

For $k \in \mathbb{Z}_+$, let $y_k = 1 - F(k) = \sum_{i > k} p_i$ be the discrete survival function ($y_m = 1$ for $m \leq -1$) and let $m(t) = \min\{j \in \mathbb{Z}_+ | y_j < 1/t\}, t \geq 0$, be the quantile function. The discrete failure rate or hazard rate $r_k$ is defined by $r_k = P[X_1 = k | X_1 \geq k] = P[X_1 = k]/P[X_1 \geq k] = p_k/y_{k-1}$, while, for $\delta \in \mathbb{Z}$, the $\delta$ failure rate is defined by

$$s_k^\delta = \frac{p_{k+\delta}}{y_{k-1}} = \frac{P[X_1 = k + \delta]}{P[X_1 \geq k]}. \tag{2.1}$$

Finally, let the cumulative $\delta$ failure rate be given by $\theta^\delta(k) = \sum_{i=0}^k s_i^\delta$ with $\theta^\delta(\infty) = \sum_{i=0}^\infty s_i^\delta \leq \infty$ and $\Theta^\delta(t) = \max\{k \in \mathbb{Z}_+ \mid \theta^\delta(k) \leq t\}$ for $t \in [s_0^\delta, \theta^\delta(\infty))$ (from now on the superscript $\delta$ is dropped for simplicity). Then $t \in [\theta(\Theta(t)), \theta(\Theta(t) + 1))$ and $P[\theta(X_n) > t] = P[X_n > \Theta(t)] = y_{\Theta(t)}$ for all $t \in [s_0, \theta(\infty))$ and $n \geq 1$.

It is easy to verify that $r_k = 1 - y_k/y_{k-1}$, $y_k = \prod_{i=0}^{k}(1-r_i)$ and, consequently,

$$s_k = r_{k+\delta}\frac{y_{k+\delta-1}}{y_{k-1}} = \begin{cases} r_{k+\delta}\prod_{i=k}^{k+\delta-1}(1-r_i), & \text{for } \delta > 0, \\ \dfrac{r_{k+\delta}}{\prod_{i=k+\delta}^{k-1}(1-r_i)}, & \text{for } \delta < 0. \end{cases} \quad (2.2)$$

Martingales are defined relative to the natural filtration of the observations $\{\mathcal{F}_n, n \geq 0\}$, with $\mathcal{F}_n = \sigma(X_1, \ldots, X_n)$ for $n \geq 1$ and $\mathcal{F}_0 = \{\varnothing, \Omega\}$. Convergence of a sequence of real numbers $\{a_n, n \geq 1\}$ to a limit $a$, as $n \to \infty$, is denoted $\lim_n a_n = a$ or $a_n \xrightarrow[n]{} a$. We write $a_n \underset{n}{\sim} b_n$ if either $a_n$ and $b_n$ both go to infinity or zero as $n \to \infty$, with $\lim_n a_n/b_n = 1$, or both converge to non-zero finite limits as $n \to \infty$. When $a_n$ diverges increasingly to infinity as $n \to \infty$, we write $a_n \uparrow \infty$. For convergence in probability and weak convergence, we use the superscripted arrows $\xrightarrow[n]{P}$ and $\xrightarrow[n]{D}$, respectively. The centered normal distribution with variance $\sigma^2$ is denoted by $N(0, \sigma^2)$.

**Proposition 2.1.** *Let $\delta \in \mathbb{Z}$, let $N_n = \sum_{k=1}^{n} I_k$ be the counting process of $\delta$-records, with $I_k = \mathbf{1}_{\{X_k > M_{k-1} + \delta\}}$, and let $\theta(k) = \sum_{i=0}^{k} s_i$, where $s_i$ is defined in (2.1). Then*

$$N_n - \theta(M_n) = N_n - \sum_{k=0}^{M_n} s_k, \qquad n \geq 1, \quad (2.3)$$

*is a martingale. Moreover, the martingale is cubic integrable if* (a) $\delta < 0$ *and* $\limsup_k r_k < 1$ *or* (b) $\delta < 0$, $\lim_k r_k = 1$ *and* $\lim_k (1-r_k)/(1-r_{k-1}) = 1$ *or* (c) $\delta \geq 0$.

**Proof.** Clearly $\mathrm{E}[I_k|\mathcal{F}_{k-1}] = P[X_k > M_{k-1} + \delta|\mathcal{F}_{k-1}] = 1 - F(M_{k-1} + \delta) = y_{M_{k-1}+\delta}$. On the other hand, letting $\Delta\theta(M_k) = \theta(M_k) - \theta(M_{k-1})$, we get

$$\mathrm{E}[\Delta\theta(M_k)|\mathcal{F}_{k-1}] = \mathrm{E}\left[\sum_{i=0}^{M_k} s_i - \sum_{i=0}^{M_{k-1}} s_i \Big| \mathcal{F}_{k-1}\right]$$

$$= \sum_{i=1}^{\infty}\sum_{j=1}^{i} s_{M_{k-1}+j} P[X_k = M_{k-1} + i|\mathcal{F}_{k-1}]$$

$$= \sum_{j=1}^{\infty} s_{M_{k-1}+j} P[X_k > M_{k-1} + j - 1|\mathcal{F}_{k-1}]$$

$$= \sum_{j=1}^{\infty} s_{M_{k-1}+j} y_{M_{k-1}+j-1} = \sum_{j=1}^{\infty} p_{M_{k-1}+j+\delta} = y_{M_{k-1}+\delta}.$$

Therefore, $N_n - \theta(M_n)$ is a martingale. For cubic integrability of (2.3), it suffices to check cubic integrability of $\theta(X_n)$:

$$E[\theta(X_n)^3] = \sum_{k=0}^{\infty} \left(\sum_{i=0}^{k} s_i\right)^3 p_k$$

$$= \sum_{k=0}^{\infty} \left(\sum_{i=0}^{k} s_i^3 + 3\sum_{i=0}^{k-1}\sum_{j=i+1}^{k} s_i^2 s_j + 3\sum_{i=0}^{k-1}\sum_{j=i+1}^{k} s_i s_j^2 + 6\sum_{i=0}^{k-2}\sum_{j=i+1}^{k-1}\sum_{l=j+1}^{k} s_i s_j s_l\right) p_k$$

$$= \sum_{i=0}^{\infty} s_i^3 \sum_{k=i}^{\infty} p_k + 3\sum_{i=0}^{\infty} s_i^2 \sum_{j=i+1}^{\infty} s_j \sum_{k=j}^{\infty} p_k + 3\sum_{i=0}^{\infty} s_i \sum_{j=i+1}^{\infty} s_j^2 \sum_{k=j}^{\infty} p_k$$

$$+ 6\sum_{i=0}^{\infty} s_i \sum_{j=i+1}^{\infty} s_j \sum_{l=j+1}^{\infty} s_l \sum_{k=l}^{\infty} p_k$$

$$= \sum_{i=0}^{\infty} s_i^3 y_{i-1} + 3\sum_{i=0}^{\infty} s_i^2 \sum_{j=i+1}^{\infty} s_j y_{j-1}$$

$$+ 3\sum_{i=0}^{\infty} s_i \sum_{j=i+1}^{\infty} s_j^2 y_{j-1} + 6\sum_{i=0}^{\infty} s_i \sum_{j=i+1}^{\infty} s_j \sum_{l=j+1}^{\infty} s_l y_{l-1}. \qquad (2.4)$$

We now show that (2.4) is finite under (a). From (2.2) and the hypothesis $\limsup_k r_k < 1$, it follows that

$$s_k \leq A r_{k+\delta} \quad \text{and} \quad y_k \leq y_{k-1} \leq B y_k \qquad (2.5)$$

for $k \in \mathbb{Z}_+$ and constants $A, B > 0$. Then (2.4) is bounded above by

$$A^3 \left(\sum_{i=0}^{\infty} r_{i+\delta} y_{i+\delta-1} + 3\sum_{i=0}^{\infty} r_{i+\delta} \sum_{j=i+1}^{\infty} r_{j+\delta} y_{j+\delta-1} + 3\sum_{i=0}^{\infty} r_{i+\delta} \sum_{j=i+1}^{\infty} r_{j+\delta} y_{j+\delta-1} \right.$$

$$\left. + 6\sum_{i=0}^{\infty} r_{i+\delta} \sum_{j=i+1}^{\infty} r_{j+\delta} \sum_{l=j+1}^{\infty} r_{l+\delta} y_{l+\delta-1}\right)$$

$$\leq A^3 \left(1 + 6\sum_{i=0}^{\infty} r_{i+\delta} y_{i+\delta-1} + 6\sum_{i=0}^{\infty} r_{i+\delta} \sum_{j=i+1}^{\infty} r_{j+\delta} y_{j+\delta-1}\right) \leq 13 A^3.$$

We now consider (b). Let $T_1 = \sum_{i=0}^{\infty} s_i^3 y_{i-1}$ and note that

$$s_i^3 y_{i-1} = \frac{p_{i+\delta}^3}{y_{i-1}^3} y_{i-1} = \frac{r_{i+\delta}^3 y_{i+\delta-1}^3}{y_{i-1}^2} \leq \frac{y_{i+\delta-1}^3}{y_{i-1}^2} \underset{i}{\sim} y_{i+3\delta-1} = \frac{p_{i+3\delta}}{r_{i+3\delta}} \underset{i}{\sim} p_{i+3\delta},$$

where the last two equivalences follow from $\lim_k (1-r_k)/(1-r_{k-1}) = 1$ and $\lim_k r_k = 1$, respectively. Hence, $T_1 < \infty$.

Let $T_2 = \sum_{i=0}^{\infty} s_i^2 \sum_{j>i}^{\infty} s_j y_{j-1}$ and note that $s_j y_{j-1} = p_{j+\delta}$ yields $T_2 = \sum_{i=0}^{\infty} s_i^2 y_{i+\delta}$. However,

$$s_i^2 y_{i+\delta} = \frac{p_{i+\delta}^2}{y_{i-1}^2} y_{i+\delta} = \frac{r_{i+\delta}^2 y_{i+\delta-1}^2 y_{i+\delta}}{y_{i-1}^2} \leq \frac{y_{i+\delta-1}^3}{y_{i-1}^2}$$

and clearly, $T_2 < \infty$.

Let $T_3 = \sum_{i=0}^{\infty} s_i \sum_{j>i}^{\infty} s_j^2 y_{j-1}$ and observe that

$$s_j^2 y_{j-1} = \frac{p_{j+\delta}^2}{y_{j-1}^2} y_{j-1} = \frac{r_{j+\delta}^2 y_{j+\delta-1}^2 y_{j-1}}{y_{j-1}^2} \leq \frac{y_{j+\delta-1}^2}{y_{j-1}} \underset{j}{\sim} y_{j+2\delta-1} \underset{j}{\sim} p_{j+2\delta}.$$

Therefore, $T_3 \leq C \sum_{i=0}^{\infty} s_i y_{i+2\delta}$ for some constant $C > 0$, but

$$s_i y_{i+2\delta} = \frac{r_{i+\delta} y_{i+\delta-1} y_{i+2\delta}}{y_{i-1}} \leq \frac{y_{i+2\delta}^2}{y_{i-1}} \underset{i}{\sim} p_{i+4\delta+2}$$

and, hence, $T_3 < \infty$.

Last, $T_4 = \sum_{i=0}^{\infty} s_i \sum_{j>i}^{\infty} s_j \sum_{l>j}^{\infty} s_l y_{l-1}$ is similarly shown to be finite, noting that $s_l y_{l-1} = p_{l+\delta}$ and

$$s_j y_{j+\delta} = \frac{r_{j+\delta} y_{j+\delta-1} y_{j+\delta}}{y_{j-1}} \leq \frac{y_{j+\delta-1}^2}{y_{j-1}} \underset{j}{\sim} p_{j+2\delta}.$$

Finally, under condition (c), note that as $\delta \geq 0$, we have $s_k \leq r_{k+\delta}$ for all $k \in \mathbb{Z}_+$,

$$\mathrm{E}[\theta(X_n)^3] = \sum_{k=0}^{\infty} \left( \sum_{i=0}^{k} s_i \right)^3 p_k \leq \sum_{k=0}^{\infty} \left( \sum_{i=0}^{k+\delta} r_i \right)^3 p_k \leq \sum_{k=0}^{\infty} \left( \delta + \sum_{i=0}^{k} r_i \right)^3 p_k$$

and therefore it suffices to show that $\sum_{k=0}^{\infty} (\sum_{i=0}^{k} r_i)^3 p_k < \infty$ or equivalently, that all terms of (2.4) are finite, with the $s_i$ replaced by the $r_i$. Indeed,

$$\sum_{i=0}^{\infty} r_i^3 y_{i-1} = \sum_{i=0}^{\infty} r_i^2 p_i \leq 1,$$

$$\sum_{i=0}^{\infty} r_i^2 \sum_{j>i}^{\infty} r_j y_{j-1} = \sum_{i=0}^{\infty} r_i^2 y_i \leq \sum_{i=0}^{\infty} r_i y_{i-1} = 1,$$

$$\sum_{i=0}^{\infty} r_i \sum_{j>i}^{\infty} r_j^2 y_{j-1} \leq \sum_{i=0}^{\infty} r_i y_i \leq 1$$

and

$$\sum_{i=0}^{\infty} r_i \sum_{j>i}^{\infty} r_j \sum_{l>j}^{\infty} r_l y_{l-1} \leq \sum_{i=0}^{\infty} r_i \sum_{j>i}^{\infty} r_j y_{j-1} \leq \sum_{i=0}^{\infty} r_i y_i \leq 1.$$

□

**Remark 2.1.** When the random variables $X_n$ have common distribution function $F$ with density $f$, it can be shown that the process

$$N_n^\delta - \int_0^{M_n} \frac{f(x+\delta)}{1-F(x)}\, dx$$

is a martingale. We believe that our methods can be applied in this case to obtain analogous limiting results.

**Proposition 2.2.** *Let $\xi_k = I_k - \Delta\theta(M_k)$, with $I_k = \mathbf{1}_{\{X_k > M_{k-1}+\delta\}}$ and $\Delta\theta(M_k) = \theta(M_k) - \theta(M_{k-1})$, $k \geq 1$. Then the increments of the process of conditional variances of martingale (2.3) are given by*

$$\mathrm{E}[\xi_k^2|\mathcal{F}_{k-1}] = \sum_{i>M_{k-1}} s_i(y_{i+\delta} + y_{i+\delta-1} - y_{i-1}) \qquad \text{for } \delta < 0$$

*and*

$$\mathrm{E}[\xi_k^2|\mathcal{F}_{k-1}] = y_{M_{k-1}+\delta}\left(1 - 2\sum_{i=1}^{\delta} s_{M_{k-1}+i}\right) + 2\sum_{i>M_{k-1}} s_i\left(y_{i+\delta} + \frac{p_{i+\delta}}{2}\right)$$
$$- 2y_{M_{k-1}+2\delta} \qquad \text{for } \delta > 0.$$

**Proof.** We have $\mathrm{E}[\xi_k^2|\mathcal{F}_{k-1}] = \mathrm{E}[I_k|\mathcal{F}_{k-1}] - 2\mathrm{E}[I_k\Delta\theta(M_k)|\mathcal{F}_{k-1}] + \mathrm{E}[\Delta\theta(M_k)^2|\mathcal{F}_{k-1}] = y_{M_{k-1}+\delta} - 2\mathrm{E}[I_k\Delta\theta(M_k)|\mathcal{F}_{k-1}] + \mathrm{E}[\Delta\theta(M_k)^2|\mathcal{F}_{k-1}]$. Writing $m$ for $M_{k-1}$, we then have

$$\mathrm{E}[\Delta\theta(M_k)^2|\mathcal{F}_{k-1}] = \mathrm{E}[((\theta(X_k) - \theta(M_{k-1}))^+)^2|\mathcal{F}_{k-1}] = \sum_{i=1}^{\infty}\left(\sum_{j=1}^{i} s_{m+j}\right)^2 p_{m+i}$$

$$= \sum_{i=1}^{\infty}\left(\sum_{j=1}^{i} s_{m+j}^2 + 2\sum_{1\leq j_1 < j_2 \leq i} s_{m+j_1}s_{m+j_2}\right) p_{m+i}$$

$$= \sum_{j=1}^{\infty} s_{m+j}^2 \sum_{i=j}^{\infty} p_{m+i} + 2\sum_{1\leq j_1<j_2<\infty} s_{m+j_1}s_{m+j_2} \sum_{i=j_2}^{\infty} p_{m+i}$$

$$= \sum_{j=1}^{\infty} s_{m+j} p_{m+j+\delta} + 2\sum_{1\leq j_1<j_2<\infty} s_{m+j_1} p_{m+j_2+\delta}$$

$$= \sum_{j>m} s_j p_{j+\delta} + 2\sum_{j>m} s_j y_{j+\delta} = 2\sum_{j>m} s_j\left(y_{j+\delta} + \frac{p_{j+\delta}}{2}\right). \qquad (2.6)$$

When $\delta < 0$, we have $I_k \Delta\theta(M_k) = \Delta\theta(M_k)$ and

$$\mathrm{E}[\xi_k^2|\mathcal{F}_{k-1}] = \mathrm{E}[\Delta\theta(M_k)^2|\mathcal{F}_{k-1}] - y_{m+\delta}$$

$$= 2\sum_{i>m} s_i\left(y_{i+\delta} + \frac{p_{i+\delta}}{2}\right) - y_{m+\delta}$$

$$= \sum_{i>m}(2s_i y_{i+\delta} + s_i p_{i+\delta} - p_{i+\delta})$$

$$= \sum_{i>m} s_i(y_{i+\delta} + y_{i+\delta-1} - y_{i-1}).$$

Otherwise, when $\delta > 0$, we obtain

$$\mathrm{E}[I_k \Delta\theta(M_k)|\mathcal{F}_{k-1}] = \sum_{j=1}^{\infty} \sum_{i=m+1}^{m+\delta+j} s_i p_{m+\delta+j}$$

$$= \sum_{i=1}^{\delta} s_{m+i} \sum_{j=1}^{\infty} p_{m+\delta+j} + \sum_{i=\delta+1}^{\infty} s_{m+i} \sum_{j=i-\delta}^{\infty} p_{m+\delta+j}$$

$$= y_{m+\delta} \sum_{i=1}^{\delta} s_{m+i} + y_{m+2\delta}$$

and, finally, $\mathrm{E}[\xi_k^2|\mathcal{F}_{k-1}] = y_{m+\delta} + 2\sum_{i>m} s_i(y_{i+\delta} + \frac{1}{2}p_{i+\delta}) - 2(y_{m+\delta}\sum_{i=1}^{\delta} s_{m+i} + y_{m+2\delta})$. □

We now give bounds on $\mathrm{E}[|\xi_k|^3|\mathcal{F}_{k-1}]$ which will be useful for checking Lyapunov's condition in the central limit theorem.

**Proposition 2.3.** *Let $\xi_k = I_k - \Delta\theta(M_k)$, $k \geq 1$. For a positive constant $C$:*

(a) *If $\delta < 0$ and $\limsup_k r_k < 1$, then $\mathrm{E}[|\xi_k|^3|\mathcal{F}_{k-1}] \leq C y_{M_{k-1}}$ for all $k \geq 1$.*
(b) *If $\delta < 0$, $\lim_k r_k = 1$ and $\lim_k (1-r_k)/(1-r_{k-1}) = 1$, then $\mathrm{E}[|\xi_k|^3|\mathcal{F}_{k-1}] \leq C y_{M_{k-1}+3\delta}$, for all $k \geq 1$.*
(c) *If $\delta > 0$, then $\mathrm{E}[|\xi_k|^3|\mathcal{F}_{k-1}] \leq C y_{M_{k-1}+\delta}$ for all $k \geq 1$.*

**Proof.** Noting that $I_k \Delta\theta(M_k) \leq \Delta\theta(M_k)$, we have

$$\mathrm{E}[|\xi_k|^3|\mathcal{F}_{k-1}] \leq \mathrm{E}[I_k|\mathcal{F}_{k-1}] + 3\mathrm{E}[\Delta\theta(M_k)|\mathcal{F}_{k-1}]$$
$$+ 3\mathrm{E}[\Delta\theta(M_k)^2|\mathcal{F}_{k-1}] + \mathrm{E}[\Delta\theta(M_k)^3|\mathcal{F}_{k-1}]$$
$$= 4\mathrm{E}[\Delta\theta(M_k)|\mathcal{F}_{k-1}] + 3\mathrm{E}[\Delta\theta(M_k)^2|\mathcal{F}_{k-1}] + \mathrm{E}[\Delta\theta(M_k)^3|\mathcal{F}_{k-1}]. \quad (2.7)$$

We first make some calculations on the terms of (2.7) which are valid for all cases (a), (b) and (c). From Proposition 2.1 and (2.6), writing $m$ for $M_{k-1}$,

$$\mathrm{E}[\Delta\theta(M_k)|\mathcal{F}_{k-1}] = y_{m+\delta}, \qquad (2.8)$$

$$\mathrm{E}[\Delta\theta(M_k)^2|\mathcal{F}_{k-1}] = 2\sum_{j>m} s_j\left(y_{j+\delta} + \frac{p_{j+\delta}}{2}\right) \leq 2\sum_{j>m} s_j(y_{j+\delta} + p_{j+\delta})$$

$$= 2 \sum_{j>m} s_j y_{j+\delta-1}. \tag{2.9}$$

For the third moment, we proceed as in Proposition 2.2, calculating $\mathrm{E}[\Delta\theta(M_k)^2|\mathcal{F}_{k-1}]$ (see also the calculations for $\mathrm{E}[\theta(X_n)^3]$ in (2.4)):

$$\mathrm{E}[\Delta\theta(M_k)^3|\mathcal{F}_{k-1}] = \mathrm{E}[((\theta(X_k) - \theta(M_{k-1}))^+)^3|\mathcal{F}_{k-1}] = \sum_{i=1}^{\infty}\left(\sum_{j=1}^{i} s_{m+j}\right)^3 p_{m+i}$$

$$= \sum_{i=1}^{\infty}\left(\sum_{j=1}^{i} s_{m+j}^3 + 3\sum_{1\leq j_1<j_2\leq i} s_{m+j_1}^2 s_{m+j_2} + 3\sum_{1\leq j_1<j_2\leq i} s_{m+j_1} s_{m+j_2}^2\right.$$

$$\left.+ 6\sum_{1\leq j_1<j_2<j_3\leq i} s_{m+j_1} s_{m+j_2} s_{m+j_3}\right) p_{m+i}$$

$$= \sum_{j>m} s_j^2 p_{j+\delta} + 3\sum_{j>m} s_j^2 y_{j+\delta} + 3\sum_{j_1>m} s_{j_1} \sum_{j_2>j_1} s_{j_2} p_{j_2+\delta}$$

$$+ 6 \sum_{j_1>m} s_{j_1} \sum_{j_2>j_1} s_{j_2} y_{j_2+\delta}$$

$$= 3\sum_{j>m} s_j^2(y_{j+\delta} + \tfrac{1}{3}p_{j+\delta}) + 6\sum_{j_1>m} s_{j_1} \sum_{j_2>j_1} s_{j_2}(y_{j_2+\delta} + \tfrac{1}{2}p_{j_2+\delta})$$

$$\leq 3\sum_{j>m} s_j^2 y_{j+\delta-1} + 6\sum_{j_1>m} s_{j_1} \sum_{j_2>j_1} s_{j_2} y_{j_2+\delta-1}. \tag{2.10}$$

Consider now (a). From (2.5) and (2.8), $\mathrm{E}[\Delta\theta(M_k)|\mathcal{F}_{k-1}] \leq B^{-\delta} y_m$. On the other hand, from (2.5) and (2.9),

$$\mathrm{E}[\Delta\theta(M_k)^2|\mathcal{F}_{k-1}] \leq 2B^{-\delta}\sum_{j>m} s_j y_{j-1} = 2B^{-\delta} y_{m+\delta} \leq 2B^{-2\delta} y_m.$$

Finally, from (2.5) and (2.10),

$$\mathrm{E}[\Delta\theta(M_k)^3|\mathcal{F}_{k-1}] \leq 3B^{-\delta}\left(\sum_{j>m} s_j^2 y_{j-1} + 2\sum_{j_1>m} s_{j_1}\sum_{j_2>j_1} s_{j_2} y_{j_2-1}\right)$$

$$= 3B^{-\delta}\left(\sum_{j>m} s_j p_{j+\delta} + 2\sum_{j>m} s_j y_{j+\delta}\right) \leq 6B^{-\delta}\sum_{j>m} s_j y_{j+\delta-1}$$

$$\leq 6B^{-2\delta}\sum_{j>m} s_j y_{j-1} = 6B^{-2\delta} y_{m+\delta} \leq 6B^{-3\delta} y_m.$$

For case (b), we have, from (2.8) and $\delta < 0$, $\mathrm{E}[\Delta\theta(M_k)|\mathcal{F}_{k-1}] \leq y_{m+3\delta}$.

From (2.2), we have

$$s_i/s_{i+\delta} = r_{i+\delta} \prod_{j=i+2\delta}^{i+\delta-1} (1-r_j) \bigg/ \left( r_{i+2\delta} \prod_{j=i+\delta}^{i-1} (1-r_j) \right) \xrightarrow[i]{} 1, \qquad (2.11)$$

so then

$$s_i y_{i+\delta-1} \underset{i}{\sim} s_{i+\delta} y_{i+\delta-1} = p_{i+2\delta}. \qquad (2.12)$$

Therefore, from (2.9) and (2.12), $\mathrm{E}[\Delta\theta(M_k)^2|\mathcal{F}_{k-1}] \leq 2\sum_{j>m} s_j y_{j+\delta-1} \leq C y_{m+2\delta} \leq C y_{m+3\delta}$. To bound (2.10), note from (2.11) and (2.12) that

$$\sum_{j>m} s_j^2 y_{j+\delta-1} \underset{m}{\sim} \sum_{j>m} s_j p_{j+2\delta} \leq \sum_{j>m} s_j y_{j+2\delta-1} \underset{m}{\sim} \sum_{j>m} s_{j+2\delta} y_{j+2\delta-1} \underset{m}{\sim} y_{m+3\delta}$$

and

$$\sum_{j_1>m} s_{j_1} \sum_{j_2>j_1} s_{j_2} y_{j_2+\delta-1} \underset{n}{\sim} \sum_{j_1>m} s_{j_1} \sum_{j_2>j_1} s_{j_2+\delta} y_{j_2+\delta-1}$$

$$\underset{m}{\sim} \sum_{j>m} s_j y_{j+2\delta} \leq \sum_{j>m} s_j y_{j+2\delta-1} \underset{m}{\sim} y_{m+3\delta}.$$

Hence, $\mathrm{E}[\Delta\theta(M_k)^3|\mathcal{F}_{k-1}] \leq C y_{m+3\delta}$.

For (c), we have to bound (2.9) and (2.10). For $\delta > 0$, we have

$$\sum_{j>m} s_j y_{j+\delta-1} \leq \sum_{j>m} s_j y_{j-1} = y_{m+\delta},$$

$$\sum_{j>m} s_j^2 y_{j+\delta-1} \leq \sum_{j>m} s_j^2 y_{j-1} = \sum_{j>m} s_j p_{j+\delta} \leq \sum_{j>m} s_j y_{j+\delta-1} \leq y_{m+\delta},$$

$$\sum_{j_1>m} s_{j_1} \sum_{j_2>j_1} s_{j_2} y_{j_2+\delta-1} \leq \sum_{j_1>m} s_{j_1} \sum_{j_2>j_1} s_{j_2} y_{j_2-1} = \sum_{j>m} s_j y_{j+\delta} \leq \sum_{j>m} s_j y_{j-1} = y_{m+\delta}. \qquad \square$$

## 3. Central limit theorems for *δ* < 0

We first show that (A.3) and (A.4) of Theorem A.2 in the Appendix hold under mild conditions on the failure rates $r_k$. We recall that $I_k = \mathbf{1}_{\{X_k > M_{k-1}+\delta\}}$ and $\Delta\theta(M_k) = \theta(M_k) - \theta(M_{k-1})$, $k \geq 1$.

**Proposition 3.1.** *Let* $\delta < 0$, $\xi_k = I_k - \Delta\theta(M_k)$ *and*

$$z_k = \sum_{i>k} s_i(y_{i+\delta} + y_{i+\delta-1} - y_{i-1}), \qquad k \geq 1. \qquad (3.1)$$

(a) *If* $\limsup_k r_k < 1$, *then* (A.3) *holds with*

$$b_n^2 = \sum_{k=0}^{m(n)} z_k r_k / y_k. \quad (3.2)$$

(b) *If* $\lim_k r_k = 1$ *and* $\lim_k (1 - r_k)/(1 - r_{k-1}) = 1$, *then* (A.3) *holds with*

$$b_n^2 = \sum_{k=0}^{m(n)} (1 - r_k)^{2\delta}. \quad (3.3)$$

**Proof.** From Proposition 2.2,

$$\mathrm{E}[\xi_k^2 | \mathcal{F}_{k-1}] = \sum_{i > M_{k-1}} s_i (y_{i+\delta} + y_{i+\delta-1} - y_{i-1}) = z_{M_{k-1}}.$$

Note that because $\delta < 0$, then $y_{i+\delta-1} \geq y_{i-1}$ and, consequently, $z_k$ is decreasing. Thus, $z_{M_{k-1}} = \min\{z_{X_1}, \ldots, z_{X_{k-1}}\}$ for $k \geq 2$, where the random variables $z_{X_k}$ are independent, identically distributed and take values $z_j$ with probabilities $p_j$. Their common distribution function is given by $G(z) = \sum_{i \geq j} p_i = y_{j-1}$ for $z_j \leq z < z_{j-1}$ and its inverse is given by $G^-(t) = z_j$ for $y_j < t \leq y_{j-1}$. Equivalently, $G^-(1/t) = z_{m(t)}$, where $m(t)$ is the quantile function defined at the beginning of Section 2.

We obtain (a) and (b) if we show

$$\sum_{k=1}^{n} z_{M_{k-1}} / b_n^2 \xrightarrow[n]{P} 1. \quad (3.4)$$

To get (3.4), we apply Deheuvels' theorem (Theorem A.1 herein). We first determine the normalizing sequence $H(\log n)$ as follows. Let $t \geq 1$. Then

$$H(\log t) = \int_1^t G^-(1/u) \, du = \int_1^t z_{m(u)} \, du = \sum_{j=0}^{m(t)} \int_{y_j^{-1}}^{y_{j-1}^{-1}} z_{m(u)} \, du - \int_t^{y_{m(t)}^{-1}} z_{m(u)} \, du$$

$$= \sum_{j=0}^{m(t)} z_j (y_j^{-1} - y_{j-1}^{-1}) - z_{m(t)} (y_{m(t)}^{-1} - t) = \sum_{j=0}^{m(t)} \frac{z_j r_j}{y_j} - \rho(t), \quad (3.5)$$

where $\rho(t) = z_{m(t)} (y_{m(t)}^{-1} - t)$.

Consider (a). From (2.5) we obtain

$$z_k \leq \sum_{i > k} s_i (y_{i+\delta} + y_{i+\delta-1}) \leq 2 \sum_{i > k} s_i y_{i+\delta-1}$$

$$\leq 2A \sum_{i > k} r_{i+\delta} y_{i+\delta-1} = 2A y_{k+\delta} \leq C y_k, \quad (3.6)$$

with $C = 2AB^{-\delta}$.

The above upper bound for $z_k$ yields immediately $\rho(t) \leq z_{m(t)}(y_{m(t)}^{-1} - y_{m(t)-1}^{-1}) = z_{m(t)}r_{m(t)}/y_{m(t)} \leq C$ and we have $H(\log n) \sim b_n^2$. It remains to check hypotheses (A.1) and (A.2) of Theorem A.1. To this end, consider the inequality

$$z_{m(t)}/y_{m(t)-1} \leq tG^-(1/t) < z_{m(t)}/y_{m(t)}, \tag{3.7}$$

which is an easy consequence of the definitions of $m(t)$ and $G^-$. On the other hand, from (2.2) and because $y_{i+\delta} \geq y_{i-1}$, it is clear that

$$z_k \geq \sum_{i>k} s_i y_{i+\delta-1} \geq \sum_{i>k} r_{i+\delta} y_{i+\delta-1} = \sum_{i>k} p_{i+\delta} = y_{k+\delta} \geq y_{k-1}. \tag{3.8}$$

Hence, from (3.6), (3.7) and (3.8),

$$1/t \leq G^-(1/t) \leq C/t \tag{3.9}$$

for all $t > 1$ and, clearly, $H(\log t)$ has a logarithmic growth to infinity as $t \to \infty$.

Finally, from the definition of $H$ and (3.9) we get

$$0 \leq (H(x_n + \log n) - H(\log n))/H(\log n) \leq Cx_n/\log n$$

for $n \geq 2$ and (A.1) follows by taking $x_n = \log(\log n + 3)$. Also, (A.2) is readily obtained from (3.9) because

$$\sum_{k=1}^n k(G^-(1/k))^2 \left(\sum_{k=1}^n G^-(1/k)\right)^{-2} \leq C^2 \sum_{k=1}^n (1/k) \left(\sum_{k=1}^n (1/k)\right)^{-2} \xrightarrow[n]{} 0.$$

Therefore, (3.4) follows from Theorem A.1.

For (b), observe that

$$z_k = \sum_{i>k} s_i y_{i+\delta-1}(y_{i+\delta}/y_{i+\delta-1} + 1 - y_{i-1}/y_{i+\delta-1}) \underset{k}{\sim} \sum_{i>k} s_i y_{i+\delta-1}$$

and, from (2.12), we have $z_k \underset{k}{\sim} \sum_{i>k} p_{i+2\delta} = y_{k+2\delta}$. Also, as in part (a),

$$\rho(t) \leq z_{m(t)} r_{m(t)}/y_{m(t)} \underset{t}{\sim} y_{m(t)+2\delta}/y_{m(t)} \underset{t}{\sim} (1 - r_{m(t)})^{2\delta}.$$

On the other hand, $\sum_{k=0}^{m(n)} z_k r_k/y_k \underset{n}{\sim} \sum_{k=0}^{m(n)} (1 - r_k)^{2\delta} \xrightarrow[n]{} \infty$, and it is clear from Lemma A.1 that $(1 - r_n)^{2\delta}/\sum_{k=0}^n (1 - r_k)^{2\delta} \xrightarrow[n]{} 0$ and, hence, $b_n^2 \underset{n}{\sim} H(\log n)$. Next we check hypothesis (A.1) of Theorem A.1, which is clearly equivalent to

$$\sum_{k=0}^{m(nu_n)} (1 - r_k)^{2\delta} \Big/ \sum_{k=0}^{m(n)} (1 - r_k)^{2\delta} \xrightarrow[n]{} 1 \tag{3.10}$$

for some sequence $u_n \uparrow \infty$. It can be shown that $m(nu_n) - m(n) - 1 < C \log u_n$ for some $C > 0$ and all $n \geq 1$. In fact, because $\lim_k r_k = 1$, there exists a constant $a > 0$ such that $1 - r_k < a < 1$ for all $k \geq 0$. Next we consider the inequalities

$$\frac{1}{nu_n} \leq y_{m(nu_n)-1} = y_{m(n)} \prod_{i=m(n)+1}^{m(nu_n)-1} (1 - r_i) < \frac{1}{n} a^{m(nu_n)-m(n)-1}$$

for all $n$ such that $m(nu_n) - m(n) \geq 1$, which implies the desired inequality. Therefore,

$$\sum_{k=m(n)+1}^{m(nu_n)} (1 - r_k)^{2\delta} \leq \sum_{k=m(n)+1}^{m(n)+\lceil C \log u_n \rceil + 1} (1 - r_k)^{2\delta}$$

and (3.10) is proved if we establish

$$\sum_{k=m(n)+1}^{m(n)+v_n} (1 - r_k)^{2\delta} \Big/ \sum_{k=0}^{m(n)} (1 - r_k)^{2\delta} \xrightarrow[n]{} 0$$

for some $v_n \uparrow \infty$ or, equivalently, because $m(n)$ is increasing,

$$\sum_{k=n+1}^{n+w_n} (1 - r_k)^{2\delta} \Big/ \sum_{k=0}^{n} (1 - r_k)^{2\delta} \xrightarrow[n]{} 0 \tag{3.11}$$

for some $w_n \uparrow \infty$.

To prove (3.11), let $c_n^{(k)} = (1 - r_{n+k})^{2\delta} / \sum_{i=0}^{n} (1 - r_i)^{2\delta}$ and for each $l \geq 1$, let $n_l$ be such that $\max\{c_n^{(k)} | k = 1, \ldots, l\} \leq 1/l^2$ for all $n \geq n_l$. This can be done for each $l$, choosing the $n_l$'s strictly increasing because $c_n^{(k)} \xrightarrow[n]{} 0$ for all $k$. We can now define the sequence $\{w_n, n \geq 1\}$ as $w_n = l$ if $n_l \leq n < n_{l+1}$.

Consider next $\varepsilon > 0$ arbitrary and choose $l$ such that $1/l < \varepsilon$. Let $n \geq n_l$. Then $n \in [n_{l+k}, n_{l+k+1})$ for some $k \geq 0$ and $w_n = l + k$, so $c_n^{(j)} < 1/(l+k)^2$ for $j = 1, \ldots, l+k$. Thus, $\sum_{j=1}^{w_n} c_n^{(j)} = \sum_{j=1}^{l+k} c_n^{(j)} < 1/(l+k) \leq 1/l < \varepsilon$ and (3.11) follows.

For condition (A.2) in Theorem A.1, note that

$$\sum_{i=1}^{n} iG^-(1/i)^2 = \sum_{k=1}^{m(n)} \sum_{i \leq n, m(i)=k} iG^-(1/i)^2 \leq \sum_{k=1}^{m(n)} \sum_{m(i)=k} iz_k^2 = \sum_{k=1}^{m(n)} z_k^2 h(k),$$

with

$$h(k) = \sum_{m(i)=k} i = \sum_{y_{k-1}^{-1} \leq i < y_k^{-1}} i = \frac{\lceil 1/y_k \rceil^2 - \lceil 1/y_{k-1} \rceil^2}{2} - \frac{\lceil 1/y_k \rceil - \lceil 1/y_{k-1} \rceil}{2}$$

$$\underset{k}{\sim} (\lceil 1/y_k \rceil^2 - \lceil 1/y_{k-1} \rceil^2)/2 \underset{k}{\sim} (y_k^{-2} - y_{k-1}^{-2})/2,$$

where the last equivalence follows from $1/y_k - 1/y_{k-1} \underset{k}{\longrightarrow} \infty$, because $\lim_k r_k = 1$. Hence,

$$\sum_{k=1}^{m(n)} z_k^2 h(k) \underset{n}{\sim} \sum_{k=1}^{m(n)} z_k^2(y_k^{-2} - y_{k-1}^{-2})/2 \underset{n}{\sim} \sum_{k=1}^{m(n)} y_{k+2\delta}^2(y_k^{-2} - y_{k-1}^{-2})/2$$

$$\underset{n}{\sim} \sum_{k=1}^{m(n)} (1-r_k)^{4\delta} r_k(2-r_k)/2 \underset{n}{\sim} \frac{1}{2}\sum_{k=1}^{m(n)} (1-r_k)^{4\delta}.$$

It is easy to see that $\sum_{k=1}^n G^-(1/k) \underset{n}{\sim} H(\log n) \underset{n}{\sim} b_n^2$ and we have, from Lemma A.1,

$$H(\log n)^{-2} \sum_{i=1}^n iG^-(1/i)^2 \leq C \sum_{k=1}^{m(n)} (1-r_k)^{4\delta} \left(\sum_{k=1}^{m(n)} (1-r_k)^{2\delta}\right)^{-2} \underset{n}{\longrightarrow} 0.$$

Hence, (3.4) follows. □

**Proposition 3.2.** *Let $\delta < 0$ and $\xi_k = I_k - \Delta\theta(M_k)$.*

  (a) *If $\limsup_k r_k < 1$, then (A.4) holds with $b_n$ defined by (3.2).*
  (b) *If $\lim_k r_k = 1$ and $\lim_k (1-r_k)/(1-r_{k-1}) = 1$, then (A.4) holds with $b_n$ defined by (3.3).*

**Proof.** (a) From Proposition 2.3(a) we have $E[|\xi_k|^3|\mathcal{F}_{k-1}] \leq Cy_{M_{k-1}}$, where $C$ is a positive constant. On the other hand, $y_{M_{k-1}} = 1 - F(M_{k-1})$ is a decreasing function of $M_{k-1}$ so that the sum in Lyapunov's condition (A.4) is bounded by $C$ times the sum of partial minima of independent identically distributed random variables taking values $y_j$ with probabilities $p_j$. Their common distribution function is denoted by $G$, with $G(y) = \sum_{i \geq j} p_i = y_{j-1}$ for $y_j \leq y < y_{j-1}$, and its inverse is denoted by $G^-(t) = y_j$ for $y_j < t \leq y_{j-1}$.

Reasoning as in Proposition 3.1(a), we obtain $\sum_{k=1}^n \min\{y_{X_1}, \ldots, y_{X_k}\}/c_n^2 \xrightarrow[n]{P} 1$, with $c_n^2 = \sum_{k=0}^{m(n)} y_k r_k / y_k = \sum_{k=0}^{m(n)} r_k$. For details, see Propositions 3.2 and 3.3 in Gouet *et al.* [9]. To conclude, note that $c_n^2 = \sum_{j=0}^{m(n)} r_j < \sum_{j=0}^{m(n)} r_j y_{j-1}/y_j \leq b_n^2$, where the second inequality comes from (3.8). Therefore, $c_n^2/b_n^3 \underset{n}{\longrightarrow} 0$ and Lyapunov's condition (A.4) follows.

(b) From Proposition 2.3(b) we have $E[|\xi_k|^3|\mathcal{F}_{k-1}] \leq Cy_{M_{k-1}+3\delta}$ and (A.4) will follow by studying the sum of partial minima

$$\sum_{k=1}^n y_{M_{k-1}+3\delta}. \tag{3.12}$$

As before, we use Theorem A.1, where calculations follow closely those in Proposition 3.1(b). We find that the scaling sequence for (3.12), denoted $\hat{b}_n$, is given by $\hat{b}_n^2 = \sum_{k=0}^{m(n)}(1-r_k)^{3\delta}$ and it can be shown, denoting $\hat{H}$ the corresponding function $H$, that $\hat{H}(\log n) \underset{n}{\sim} \sum_{k=0}^{m(n)} r_k y_{k+3\delta}/y_k \underset{n}{\sim} \hat{b}_n^2$. Conditions (A.1) and (A.2) are analogously checked and we conclude that $\sum_{k=0}^{n} \mathrm{E}[\Delta\theta(M_k)^3|\mathcal{F}_{k-1}]/\hat{b}_n^2 \xrightarrow[n]{P} 1$. Lyapunov's condition follows if $\hat{b}_n^2/b_n^3 \xrightarrow[n]{} 0$ or, equivalently, if

$$\left(\sum_{k=0}^{n}(1-r_k)^{3\delta}\right)^2 \left(\sum_{k=0}^{n}(1-r_k)^{2\delta}\right)^{-3} \xrightarrow[n]{} 0,$$

but this convergence follows from Cauchy–Schwarz inequality and Lemma A.1 because

$$\frac{(\sum_{k=0}^{n}(1-r_k)^{3\delta})^2}{(\sum_{k=0}^{n}(1-r_k)^{2\delta})^3} \leq \frac{\sum_{k=0}^{n}(1-r_k)^{2\delta}}{\sum_{k=0}^{n}(1-r_k)^{2\delta}} \frac{\sum_{k=0}^{n}(1-r_k)^{4\delta}}{(\sum_{k=0}^{n}(1-r_k)^{2\delta})^2} \xrightarrow[n]{} 0. \qquad \square$$

We now state and prove the central limit theorem for $\delta < 0$.

**Theorem 3.1.** *Let $\delta < 0$ and let $z_k$ be as defined in* (3.1).

(a) *If $\limsup_k r_k < 1$, then*

$$\frac{N_n - \theta(m(n))}{\sqrt{\sum_{k=0}^{m(n)} z_k r_k/y_k}} \xrightarrow[n]{D} N(0,1). \tag{3.13}$$

(b) *If $\lim_k r_k = 1$ and $\lim_k (1-r_k)/(1-r_{k-1}) = 1$, then*

$$\frac{N_n - \theta(m(n))}{\sqrt{\sum_{k=0}^{m(n)}(1-r_k)^{2\delta}}} \xrightarrow[n]{D} N(0,1). \tag{3.14}$$

**Proof.** (a) Using results in Propositions 3.1(a) and 3.2(a) and Theorem A.2, we have $(N_n - \theta(M_n))/b_n \xrightarrow[n]{D} N(0,1)$, with $b_n$ defined in (3.2), so (3.13) follows if we show

$$(\theta(M_n) - \theta(m(n)))/b_n \xrightarrow[n]{P} 0. \tag{3.15}$$

This will be done by comparison with the analogous result for usual records ($\delta = 0$) contained in Proposition 3 of Gouet *et al.* [10]. From (2.5) we get

$$|\theta(M_n) - \theta(m(n))| = \sum_{i=(M_n \wedge m(n))+1}^{M_n \vee m(n)} s_i \leq A \sum_{i=(M_n \wedge m(n))+1}^{M_n \vee m(n)} r_{i+\delta}$$

$$\leq A \sum_{i=(M_n \wedge m(n))+1+\delta}^{M_n \vee m(n)} r_i \leq A \sum_{i=(M_n \wedge m(n))+1}^{M_n \vee m(n)} r_i - A\delta. \tag{3.16}$$

Let $\theta^0(k) = \sum_{i=0}^{k} r_i$ be the centering function $\theta$ of the martingale for 0 records and let $b_{0n}$ be the corresponding scaling sequence defined by (3.6) in Gouet *et al.* [10]. Then, from (3.16), $|\theta(M_n) - \theta(m(n))| \leq A(|\theta^0(M_n) - \theta^0(m(n))| - \delta)$. In Propositions 2 and 3 of Gouet *et al.* [10] we find, respectively, that $b_{0n}^2$ has logarithmic growth and that $(\theta^0(M_n) - \theta^0(m(n)))/b_{0n} \xrightarrow[n]{P} 0$. Now, it is clear that (3.15) follows because, by (3.9), $b_n^2$ has logarithmic growth as well.

(b) From Propositions 3.1(b) and 3.2(b) and Theorem A.2, we obtain

$$(N_n - \theta(M_n))/b_n \xrightarrow[n]{D} N(0,1),$$

where $b_n$ is defined in (3.3). The result will follow if we show that

$$(\theta(M_n) - \theta(m(n)))/b_n \xrightarrow[n]{P} 0.$$

To that end, define $c_n^2 = \sum_{k=0}^{m(n)} s_k^2$ and note that $b_n \underset{n}{\sim} c_n$. Therefore, according to Corollary A.1, we have to establish

$$ny_{\Theta(\varepsilon c_n + \theta(m(n)))} \xrightarrow[n]{} 0 \quad \text{and} \quad ny_{\Theta(-\varepsilon c_n + \theta(m(n)))} \xrightarrow[n]{} \infty$$

for every $\varepsilon > 0$. Let then $\varepsilon > 0$. Noting that $s_{k+1}/s_k \xrightarrow[k]{} 1$, from Lemma A.1 we have $s_{m(n)+1}^2 / \sum_{k=0}^{m(n)} s_k^2 \xrightarrow[n]{} 0$ and this implies the existence of $N \in \mathbb{N}$ such that $\varepsilon^2 \sum_{k=0}^{m(n)} s_k^2 \geq s_{m(n)+1}^2$ for all $n > N$. Therefore, because $\Theta$ is increasing, we obtain

$$\Theta(\varepsilon c_n + \theta(m(n))) = \Theta\left(\varepsilon \left(\sum_{k=0}^{m(n)} s_k^2\right)^{1/2} + \sum_{k=0}^{m(n)} s_k\right) \geq \Theta\left(\sum_{k=0}^{m(n)+1} s_k\right) = m(n) + 1. \quad (3.17)$$

Moreover, using Lemma A.1 it is also possible to find $N' \in \mathbb{N}$ such that $\varepsilon^2 \sum_{k=0}^{m(n)} s_k^2 > 4(s_{m(n)} \vee s_{m(n)-1})^2$ for all $n > N'$, implying $\varepsilon^2 \sum_{k=0}^{m(n)} s_k^2 > (s_{m(n)} + s_{m(n)-1})^2$. It follows from the previous inequality that

$$\Theta(-\varepsilon c_n + \theta(m(n))) = \Theta\left(-\varepsilon \left(\sum_{k=0}^{m(n)} s_k^2\right)^{1/2} + \sum_{k=0}^{m(n)} s_k\right) \leq \Theta\left(\sum_{k=0}^{m(n)-2} s_k\right) = m(n) - 2. \quad (3.18)$$

From (3.17) and (3.18), and recalling that $y_{m(n)} < 1/n \leq y_{m(n)-1}$, we obtain

$$ny_{\Theta(\varepsilon c_n + \theta(m(n)))} \leq ny_{m(n)+1} < 1 - r_{m(n)+1} \xrightarrow[n]{} 0,$$

$$ny_{\Theta(-\varepsilon c_n + \theta(m(n)))} \geq ny_{m(n)-2} \geq \frac{ny_{m(n)-1}}{1 - r_{m(n)-1}} \geq \frac{1}{1 - r_{m(n)-1}} \xrightarrow[n]{} \infty,$$

and (3.14) is proved. □

The case of converging failure rates is detailed in the following corollary.

**Corollary 3.1.** *Let $\delta < 0$. If $\lim_k r_k = r \in [0,1)$, then*

$$(\log n)^{-1/2}(N_n - \theta(m(n))) \xrightarrow[n]{D} N(0, \sigma_r^2),$$

*where $\sigma_r^2 = -r(1-r)^\delta((1-r)^{\delta+1} + (1-r)^\delta - 1)/\log(1-r)$ if $r \neq 0$ and $\sigma_0 = 1$. Moreover:*

(a) *If $r > 0$ and $\sum_{i=0}^n |r_i - r|/\sqrt{n} \xrightarrow[n]{} 0$, then*

$$(\log n)^{-1/2}(N_n + r(1-r)^\delta \log n/\log(1-r)) \xrightarrow[n]{D} N(0, \sigma_r^2).$$

(b) *If $r = 0$ and $\sum_{i=0}^n r_i^2/\sqrt{n} \xrightarrow[n]{} 0$, then*

$$(\log n)^{-1/2}(N_n - \log n) \xrightarrow[n]{D} N(0,1).$$

**Proof.** Let us show that $b_n^2/\log n \xrightarrow[n]{} \sigma_r^2$. First note that, from identity $y_k = \prod_{i=0}^k (1-r_i)$ and the definition of $m(t)$, we have $y_{m(n)} < 1/n \leq y_{m(n)-1}$ and

$$-\sum_{k=0}^{m(n)-1} \log(1 - r_k) \leq \log n < -\sum_{k=0}^{m(n)} \log(1 - r_k). \tag{3.19}$$

For $r \in (0,1)$, let $L = (1-r)^\delta((1-r)^{\delta+1} + (1-r)^\delta - 1)$. We study the asymptotic behaviour of the three sums in the definition of $z_k$ in (3.1), for $\lim_k r_k = r$. For the first sum we obtain $\sum_{i>k} s_i y_{i+\delta} \underset{k}{\sim} (1-r)^{\delta+1} \sum_{i>k} s_i y_{i-1} = (1-r)^{\delta+1} y_{k+\delta} \underset{k}{\sim} (1-r)^{2\delta+1} y_k$. For the next, we get $\sum_{i>k} s_i y_{i+\delta-1} \underset{k}{\sim} (1-r)^{2\delta} y_k$, and for the last, $\sum_{i>k} s_i y_{i-1} = y_{k+\delta} \underset{k}{\sim} (1-r)^\delta y_k$. Collecting the above results, we find that $z_k/y_k \xrightarrow[k]{} L$ and $b_n^2 \underset{n}{\sim} L \sum_{j=0}^{m(n)} r_j \underset{n}{\sim} rLm(n)$. Finally, dividing (3.19) by $m(n)$ and taking limits, we get $\log n/m(n) \xrightarrow[n]{} -\log(1-r)$ and the conclusion follows.

Consider now the case $r = 0$. Clearly $z_k = \sum_{i>k} s_i(y_{i+\delta} + y_{i+\delta-1} - y_{i-1}) \underset{k}{\sim} \sum_{i>k} s_i y_{i-1} = y_{k+\delta} \underset{k}{\sim} y_k$. Then $b_n^2 = \sum_{j=0}^{m(n)} z_j r_j/y_j \underset{n}{\sim} \sum_{j=0}^{m(n)} r_j$. Therefore, by (3.19), $b_n^2 \underset{n}{\sim} \sum_{k=0}^{m(n)} r_k \underset{n}{\sim} -\sum_{k=0}^{m(n)} \log(1 - r_k) \underset{n}{\sim} \log n$.

We now prove (a) and (b) about the simplification of the centering sequences.

(a) When $0 < r < 1$, we have to show

$$(\log n)^{-1/2}(\theta(m(n)) + r(1-r)^\delta \log n/\log(1-r)) \xrightarrow[n]{} 0.$$

From (3.19), we have $m(n) \underset{n}{\sim} -\log n/\log(1-r)$. On the other hand, from the definition of $m(n)$, we get $y_{m(n)} < 1/n \le y_{m(n)-1}$ and

$$-1 - \frac{\sum_{i=0}^{m(n)} \log(1-R_i)}{\log(1-r)} < m(n) + \frac{\log n}{\log(1-r)} \le \frac{-\sum_{i=0}^{m(n)} \log(1-R_i)}{\log(1-r)} + \frac{\log(1-R_{m(n)})}{\log(1-r)},$$

where $R_i = (r_i - r)/(1-r)$. Dividing by $\sqrt{m(n)}$, we find that the left and right terms above tend to 0 as $n \to \infty$, obtaining thus

$$(\log n)^{-1/2}(m(n) + \log n/\log(1-r)) \underset{n}{\longrightarrow} 0. \qquad (3.20)$$

Finally, it remains to check that $(\log n)^{-1/2}(\theta(m(n)) - r(1-r)^\delta m(n)) \underset{n}{\longrightarrow} 0$, or, equivalently, $\sum_{k=0}^{n}(s_k - r(1-r)^\delta)/\sqrt{n} \underset{n}{\longrightarrow} 0$. This last convergence is obtained from an inductive argument on $-\delta$ as follows (we write the superscript $\delta$ on $s_k$ to avoid confusion). Recalling that $s_k^\delta = r_{k+\delta}/\prod_{i=k+\delta}^{k-1}(1-r_i)$ for $\delta < 0$, define $D_k^{(\delta)} = s_k^\delta - r(1-r)^\delta$. Then, for $\delta = -1$, we have

$$D_k^{(-1)} = \frac{r_{k-1}}{1 - r_{k-1}} - \frac{r}{1-r} = \frac{r_{k-1} - r}{(1-r)(1-r_{k-1})},$$

which, together with the hypothesis on the $r_k$'s, implies $\sum_{k=0}^{n} D_k^{(-1)}/\sqrt{n} \underset{n}{\longrightarrow} 0$.

Let us assume now that convergence holds for $\delta \in \mathbb{Z}_-$ and consider $D_k^{(\delta-1)}$. It is easy to see that $D_k^{(\delta-1)} = s_{k-1}^\delta/(1-r_{k-1}) - r(1-r)^\delta/(1-r)$, which, after some algebraic manipulation, yields

$$D_k^{(\delta-1)} = \frac{(1-r)D_{k-1}^{(\delta)} + r(1-r)^\delta(r_{k-1} - r)}{(1-r)(1-r_{k-1})} = \frac{D_{k-1}^{(\delta)}}{1 - r_{k-1}} + r(1-r)^\delta D_k^{(-1)}. \qquad (3.21)$$

From the inductive hypothesis and (3.21), we finally obtain $\sum_{k=0}^{n} D_k^{(\delta-1)}/\sqrt{n} \underset{n}{\longrightarrow} 0$.

(b) For $r = 0$, we have to show $|\theta(m(n)) - \log n|/\sqrt{\log n} \underset{n}{\longrightarrow} 0$, provided that $\sum_{i=1}^{n} r_i^2/\sqrt{n} \underset{n}{\longrightarrow} 0$. To that end we write

$$|\theta(m(n)) - \log n| \le \left|\theta(m(n)) - \sum_{k=0}^{m(n)} r_{k+\delta}\right| + \left|\sum_{k=0}^{m(n)} r_k - \sum_{k=0}^{m(n)} r_{k+\delta}\right| + \left|\log n - \sum_{k=0}^{m(n)} r_k\right| \qquad (3.22)$$

and show that all terms on the right of (3.22) divided by $\sqrt{m(n)}$ tend to 0 as $n \to \infty$.

Note that $|\log y_n + \sum_{k=0}^{n} r_k| \le C \sum_{k=0}^{n} r_k^2$. Then $|\log n - \sum_{k=0}^{m(n)} r_k|/\sqrt{m(n)} \underset{n}{\longrightarrow} 0$. For the second term we have $\sum_{k=0}^{n} r_k - \sum_{k=0}^{n} r_{k+\delta} = \sum_{k=n+\delta+1}^{n} r_k \le -\delta$. Furthermore, for the first term we use an inductive reasoning as done for $r > 0$ above. Let $D_k^{(\delta)} = s_k^\delta -$

$r_{k+\delta}$. Then, for $\delta = -1$, $D_k^{(-1)} = r_{k-1}/(1 - r_{k-1}) - r_{k-1} = r_{k-1}^2/(1 - r_{k-1})$ and clearly $\sum_{k=0}^n |D_k^{(-1)}|/\sqrt{n} \xrightarrow[n]{} 0$. Let us assume $\sum_{k=0}^n |D_k^{(\delta)}|/\sqrt{n} \xrightarrow[n]{} 0$. Then

$$|D_k^{(\delta-1)}| \leq |D_{k-1}^{(\delta)}|/(1 - r_{k-1}) + r_{k-1} r_{k+\delta-1}/(1 - r_{k-1}). \tag{3.23}$$

The Cauchy–Schwarz inequality applied to (the sum over $k$ of) the last term of (3.23) and the inductive hypothesis yields, finally, $\sum_{k=0}^n |D_k^{(\delta-1)}|/\sqrt{n} \xrightarrow[n]{} 0$. □

## 4. Central limit theorems for $\delta > 0$

In the following two propositions we check conditions (A.3) and (A.4) of the martingale central limit theorem for positive $\delta$-records. Attention is restricted to converging failure rates $r_k$ to reduce the study of conditional variances to sums of minima. We recall again that $I_k = \mathbf{1}_{\{X_k > M_{k-1}+\delta\}}$ and $\Delta \theta(M_k) = \theta(M_k) - \theta(M_{k-1})$, $k \geq 1$.

**Proposition 4.1.** *Let $\delta > 0$, $\lim_k r_k = r \in [0,1]$ and $\xi_k = I_k - \Delta \theta(M_k)$.*

(a) *If $r < 1$, (A.3) holds with $b_n^2 = \sigma_r^2 \log n$, where $\sigma_0 = 1$ and*

$$\sigma_r^2 = -r(1-r)^\delta[(1-r)^{\delta+1} - (1 + 2\delta r)(1-r)^\delta + 1]/\log(1-r), \text{ for } r \neq 0.$$

(b) *If $r = 1$ and $\sum_{k=1}^\infty e_k = \infty$, with $e_k = (1 - r_k) \cdots (1 - r_{k+\delta-1})$, (A.3) holds with $b_n^2 = \sum_{k=1}^{m(n)} e_k$. When $\sum_{k=1}^\infty e_k < \infty$, $\lim_n N_n < \infty$ almost surely.*

**Proof.** (a) From Proposition 2.2,

$$\mathrm{E}[\xi_k^2 | \mathcal{F}_{k-1}] = \sum_{i > M_{k-1}} (s_i(y_{i+\delta} + y_{i+\delta-1} + y_{i-1}) - 2p_{i+2\delta}) - 2y_{M_{k-1}+\delta} \sum_{i=M_{k-1}+1}^{M_{k-1}+\delta} s_i. \tag{4.1}$$

We first show that

$$\frac{\sum_{i>m}(s_i(y_{i+\delta} + y_{i+\delta-1} + y_{i-1})) - 2\sum_{i>m} p_{i+2\delta} - 2y_{m+\delta} \sum_{i=m+1}^{m+\delta} s_i}{y_m} \xrightarrow[m]{} L, \tag{4.2}$$

where $L = (1-r)^\delta((1-r)^{\delta+1} - (1+2\delta r)(1-r)^\delta + 1)$ for $r > 0$ and $L = 1$ for $r = 0$. Note that as $y_i/y_{i-1} \xrightarrow[i]{} 1 - r$, we have

$$s_i(y_{i+\delta} + y_{i+\delta-1} + y_{i-1}) = s_i y_{i-1}(y_{i+\delta}/y_{i-1} + y_{i+\delta-1}/y_{i-1} + 1)$$
$$\underset{i}{\sim} s_i y_{i-1}((1-r)^{\delta+1} + (1-r)^\delta + 1)$$
$$= p_{i+\delta}((1-r)^{\delta+1} + (1-r)^\delta + 1).$$

Then

$$\sum_{i>m} s_i(y_{i+\delta} + y_{i+\delta-1} + y_{i-1}) \underset{m}{\sim} y_{m+\delta}((1-r)^{\delta+1} + (1-r)^\delta + 1)$$

$$\underset{m}{\sim} y_m(1-r)^\delta((1-r)^{\delta+1} + (1-r)^\delta + 1).$$

Also, $\sum_{i>m} p_{i+2\delta} = y_{m+2\delta} \underset{m}{\sim} (1-r)^{2\delta} y_m$. Finally,

$$y_{m+\delta} \sum_{i=m+1}^{m+\delta} s_i \underset{m}{\sim} (1-r)^\delta y_m \sum_{i=m+1}^{m+\delta} r_{i+\delta} \prod_{k=i}^{i+\delta-1} (1-r_k) \underset{m}{\sim} (1-r)^{2\delta} \delta r y_m$$

and (4.2) is proved.

On the other hand, by Propositions 3.2 and 3.3 of Gouet *et al.* [9], we have

$$\sum_{k=1}^n y_{M_{k-1}} / \log n \xrightarrow{P}_n -r/\log(1-r) \tag{4.3}$$

for $r \in [0, 1)$, with $-0/\log 1 = 1$, and (a) is proved.

(b) Recalling expression (4.1), we first show that

$$\frac{\sum_{i>m}(s_i(y_{i+\delta} + y_{i+\delta-1} + y_{i-1})) - 2\sum_{i>m} p_{i+2\delta} - 2y_{m+\delta} \sum_{i=m+1}^{m+\delta} s_i}{y_{m+\delta}} \xrightarrow{m} 1. \tag{4.4}$$

Note that

$$\sum_{i>m} s_i(y_{i+\delta} + y_{i+\delta-1} + y_{i-1}) = \sum_{i>m} s_i y_{i-1}\left(\frac{y_{i+\delta}}{y_{i-1}} + \frac{y_{i+\delta-1}}{y_{i-1}} + 1\right)$$

$$\underset{m}{\sim} \sum_{i>m} s_i y_{i-1} = y_{m+\delta},$$

$$\sum_{i>m} \frac{p_{i+2\delta}}{y_{m+\delta}} = \frac{y_{m+2\delta}}{y_{m+\delta}} \xrightarrow{m} 0 \quad \text{and} \quad y_{m+\delta} \sum_{i=m+1}^{m+\delta} \frac{s_i}{y_{m+\delta}} = \sum_{i=m+1}^{m+\delta} r_{i+\delta} e_i \xrightarrow{m} 0.$$

Then (4.4) is proved.

Therefore, $\sum_{k=1}^n E[\xi_k^2|\mathcal{F}_{k-1}] \underset{n}{\sim} \sum_{k=1}^n y_{M_{k-1}+\delta}$ almost surely. Define the decreasing sequence $z_k = y_{k+\delta}$, $k \geq 1$. Then $\sum_{k=1}^n y_{M_{k-1}+\delta} = \sum_{k=1}^n \min\{z_{X_1}, \ldots, z_{X_k}\}$, where the random variables $z_{X_k}$ are independent, identically distributed and take values $z_j$ with probabilities $p_j$. Their common distribution function is $G(z) = \sum_{i \geq j} p_i = y_{j-1}$, $z_j \leq z < z_{j-1}$, and its inverse $G^-(t) = z_j$, $y_j < t \leq y_{j-1}$. We now apply Theorem A.1 to the sum of minima. From (3.5), we have $H(\log t) = \sum_{j=0}^{m(t)} z_j r_j / y_j - \rho(t)$. In this case,

$$\sum_{j=0}^{m(t)} \frac{z_j r_j}{y_j} \underset{t}{\sim} \sum_{j=0}^{m(t)} e_{j+1}$$

and $|\rho(t)| \leq y_{m(t)+\delta}/y_{m(t)} \leq 1$, so $H(\log n) \underset{n}{\sim} b_n^2$. Then, from Theorem A.1, if $\sum_{n=1}^{\infty} e_n < \infty$, we have $\sum_{k=1}^{\infty} y_{M_{k-1}+\delta} < \infty$ almost surely. Thus, $\sum_{k=1}^{\infty} \mathrm{E}[I_k|\mathcal{F}_{k-1}] = \sum_{k=1}^{\infty} y_{M_{k-1}+\delta} < \infty$ and from the conditional Borel–Cantelli lemma (see Neveu [17], Corollary VII-2-6), we conclude that $\lim_n N_n < \infty$.

Let now $\sum_{n=1}^{\infty} e_n = \infty$. We check hypotheses (A.1) and (A.2) of Theorem A.1. As in the proof of (A.1) in Proposition 3.1(b), it suffices to show in this case that

$$\sum_{i=m(n)+1}^{m(nv_n)} e_i \bigg/ \sum_{i=1}^{m(n)} e_i \underset{n}{\longrightarrow} 0 \qquad (4.5)$$

for some $v_n \uparrow \infty$. Because $e_n < 1$ and $m(nv_n) - m(n) - 1 < C \log v_n$ for some $C > 0$ and every $n \geq 1$, (4.5) holds taking $v_n = \sum_{i=1}^{m(n)} e_i$.

We now study (A.2) and again, as in the proof of Proposition 3.1(b), we have

$$\sum_{k=1}^{m(n)} z_k^2 h(k) \underset{n}{\sim} \frac{1}{2} \sum_{k=1}^{m(n)} y_{k+\delta}^2 (y_k^{-2} - y_{k-1}^{-2}) \underset{n}{\sim} \frac{1}{2} \sum_{k=1}^{m(n)} e_{k+1}^2 \underset{n}{\sim} \frac{1}{2} \sum_{k=1}^{m(n)} e_k^2.$$

Therefore, because $e_k < 1$,

$$H(\log n)^{-2} \sum_{i=1}^{n} i G^-(1/i)^2 \leq C \sum_{k=1}^{m(n)} e_k^2 \bigg/ \bigg(\sum_{k=1}^{m(n)} e_k\bigg)^2 \underset{n}{\longrightarrow} 0.$$

Hence, (A.3) holds because

$$\sum_{k=1}^{n} y_{M_{k-1}+\delta} \bigg/ \sum_{k=1}^{m(n)} e_k \xrightarrow[n]{P} 1. \qquad (4.6)$$

$\square$

**Proposition 4.2.** *Let $\delta > 0$, $\lim_k r_k = r \in [0,1]$ and $\xi_k = I_k - \Delta\theta(M_k)$.*

(a) *If $r < 1$, then (A.4) holds with $b_n^2 = \log n$.*
(b) *If $r = 1$ and $\sum_{k=1}^{\infty} e_k = \infty$, with $e_k = (1-r_k)\cdots(1-r_{k+\delta-1})$, then (A.4) holds with $b_n^2 = \sum_{k=1}^{m(n)} e_k$.*

**Proof.** (a) From Proposition 2.3, we have $\mathrm{E}[|\xi_k|^3|\mathcal{F}_{k-1}] \leq C y_{M_{k-1}+\delta} \leq C y_{M_{k-1}}$ for some $C > 0$. From (4.3), $\sum_{k=1}^{n} y_{M_{k-1}}$ has logarithmic growth and (A.4) holds.

(b) From Proposition 2.3, $\mathrm{E}[|\xi_k|^3|\mathcal{F}_{k-1}] \leq C y_{M_{k-1}+\delta}$ for some $C > 0$. Then (A.4) holds by (4.6). □

We now state and prove the central limit theorem for $\delta > 0$.

**Theorem 4.1.** *Let $\delta > 0$ and $\lim_k r_k = r \in [0, 1]$.*

(a) *If $r < 1$, then*

$$(\log n)^{-1/2}(N_n - \theta(m(n))) \xrightarrow[n]{D} N(0, \sigma_r^2),$$

*where $\sigma_r^2 = -r(1-r)^\delta((1-r)^{\delta+1} - (1+2\delta r)(1-r)^\delta + 1)/\log(1-r)$ for $r \neq 0$ and $\sigma_0 = 1$.*

(b) *If $r = 1$, then, defining $e_k = (1 - r_k) \cdots (1 - r_{k+\delta-1})$, we have*

$$\frac{N_n - \theta(m(n))}{\sqrt{\sum_{k=0}^{m(n)} e_k}} \xrightarrow[n]{D} N(0, 1)$$

*whenever $\sum_{k=0}^\infty e_k = \infty$ and $\lim_n N_n < \infty$ almost surely when $\sum_{k=0}^\infty e_k < \infty$.*

**Proof.** (a) By Propositions 4.1(a) and 4.2(a) and Theorem A.2, it only remains to show

$$(\log n)^{-1/2}(\theta(M_n) - \theta(m(n))) \xrightarrow[n]{P} 0. \tag{4.7}$$

We have

$$|\theta(M_n) - \theta(m(n))| = \sum_{i=(M_n \wedge m(n))+1}^{M_n \vee m(n)} s_i \leq \sum_{i=(M_n \wedge m(n))+1}^{M_n \vee m(n)} r_{i+\delta} \leq \delta + \sum_{i=(M_n \wedge m(n))+1}^{M_n \vee m(n)} r_i,$$

so $|\theta(M_n) - \theta(m(n))| \leq |\theta^0(M_n) - \theta^0(m(n))| + \delta$, with $\theta^0(k) = \sum_{i=0}^k r_i$ and, as in the proof of Theorem 3.1(a), $(\theta^0(M_n) - \theta^0(m(n)))/\sqrt{\log n} \xrightarrow[n]{P} 0$. Then (4.7) holds.

(b) By Propositions 4.1(b) and 4.2(b) and Theorem A.2, we have to prove that $(\theta(M_n) - \theta(m(n)))/(\sum_{k=1}^n e_k)^{1/2} \xrightarrow[n]{P} 0$ when $\sum_{n=1}^\infty e_n = \infty$. This follows from inequality $|\theta(M_n) - \theta(m(n))| \leq |M_n - m(n)|$ and the tightness of $M_n - m(n)$ when $\lim_k r_k \to 1$ (see the proof of Theorem 1 in Gouet *et al.* [10]). □

**Corollary 4.1.** *Under the hypotheses of Theorem 4.1(a) we have*

(a) *If $r > 0$ and $\sum_{i=0}^n |r_i - r|/\sqrt{n} \xrightarrow[n]{} 0$, then*

$$(\log n)^{-1/2}\left(N_n + \frac{r(1-r)^\delta \log n}{\log(1-r)}\right) \xrightarrow[n]{D} N(0, \sigma_r^2).$$

(b) If $r=0$ and $\sum_{i=0}^{n} r_i^2/\sqrt{n} \xrightarrow[n]{} 0$, then

$$(\log n)^{-1/2}(N_n - \log n) \xrightarrow[n]{D} N(0,1).$$

**Proof.** The proof is very similar to the proof of Corollary 3.1 except for some changes in our inductive arguments.

(a) Recalling that $s_k^\delta = r_{k+\delta} \prod_{i=k}^{k+\delta-1}(1-r_i)$, define $D_k^{(\delta)} = s_k^\delta - r(1-r)^\delta$. From (3.20), we have to prove that $\sum_{i=0}^{n}|r_i - r|/\sqrt{n} \xrightarrow[n]{} 0$ implies $\sum_{k=0}^{n} D_k^{(\delta)}/\sqrt{n} \xrightarrow[n]{} 0$ for $\delta = 1, 2, \ldots$. For $\delta = 1$, we have

$$D_k^{(1)} = r_{k+1} - r + r^2 - r_k r_{k+1} = r_{k+1} - r + r(r - r_k) + r_k(r - r_{k+1})$$

and $\sum_{k=0}^{n} D_k^{(1)}/\sqrt{n} \xrightarrow[n]{} 0$.

Assume $\sum_{k=0}^{n} D_k^{(\delta)}/\sqrt{n} \xrightarrow[n]{} 0$ and note that

$$D_k^{(\delta+1)} = s_k^{\delta+1} - r(1-r)^{\delta+1} = r_{k+\delta+1} \prod_{i=k}^{k+\delta}(1-r_i) - r(1-r)^{\delta+1}$$

$$= s_{k+1}^\delta(1-r_k) - r(1-r)^{\delta+1} = (1-r_k)D_{k+1}^{(\delta)} + r(1-r)^\delta(r - r_k).$$

Then, clearly, $\sum_{k=0}^{n} D_k^{(\delta+1)}/\sqrt{n} \xrightarrow[n]{} 0$.

(b) We prove that $\sum_{k=1}^{n}|D_k^{(\delta)}|/\sqrt{n} \xrightarrow[n]{} 0$ under $\sum_{i=1}^{n} r_i^2/\sqrt{n} \xrightarrow[n]{} 0$, where $D_k^{(\delta)} = s_k^\delta - r_{k+\delta}$. For $\delta = 1$, we have $D_k^{(1)} = r_{k+1}(1-r_k) - r_{k+1} = -r_k r_{k+1}$, and $\sum_{k=1}^{n}|D_k^{(1)}|/\sqrt{n} \xrightarrow[n]{} 0$ follows from the Cauchy–Schwarz inequality. Consider the inductive hypothesis $\sum_{k=1}^{n}|D_k^{(\delta)}|/\sqrt{n} \xrightarrow[n]{} 0$. Then $D_k^{(\delta+1)} = s_k^{\delta+1} - r_{k+\delta+1} = s_{k+1}^\delta(1-r_k) - r_{k+\delta+1} = D_{k+1}^{(\delta)} - r_k s_{k+1}^\delta$ and

$$\sum_{k=1}^{n}|D_k^{(\delta+1)}| \le \sum_{k=1}^{n}|D_k^{(\delta)}| + \sum_{k=1}^{n} r_k s_{k+1}^\delta,$$

which tends to 0 divided by $\sqrt{n}$ from the inductive hypothesis and the Cauchy–Schwarz inequality, because $s_{k+1}^\delta \le r_{k+1+\delta}$. $\square$

*Remark 4.1.* Notice that Theorem 4.1(a) is more restrictive than Theorem 3.1(a), concerning the behaviour of the failure rates $r_k$. This is because the process of conditional variances (A.3) can be written as partial sums of minima only when $\delta < 0$ (see Proposition 2.2). For positive $\delta$, we were able to analyze the case of converging $r_k$'s, where conditional variances behave asymptotically as sums of minima.

On the other hand, comparing Theorem 3.1(b) and Theorem 4.1(b) about distributions with light tails ($\lim_k r_k = 1$), we find more generality in the positive case because we do not impose any condition on the rate of convergence of $r_k$ to 1. This is not surprising in view of the structure of the $\delta$ failure rates $s_k$, with $1 - r_k$'s in the denominator when $\delta$ is negative. In this case, it can be shown that, for the martingale central limit theorem, Theorem A.2, it is enough to have $(1 - r_k)/(1 - r_{k-1})$ bounded away from zero and infinity; however, the change of the centering sequence $\theta(M_n)$ by a deterministic one needs some extra hypothesis on the convergence of $r_k$ to 1.

**Remark 4.2.** When $\delta > 0$, unlike the negative case, it is not guaranteed that the number of $\delta$-records is infinite. Nevertheless, when this happens, this number is always asymptotically normal in contrast to the situation of usual records, which can grow to infinity without having an asymptotically normal distribution; see Gouet *et al.* ([10], Theorem 1(b)).

## 5. Examples

**Example 5.1** *(Geometric)*. We consider independent identically distributed random variables with geometric distribution on $\mathbb{Z}_+$, that is, $p_k = pq^k, k \geq 0, n \geq 1$, with $p \in (0,1)$ and $q = 1 - p$. Clearly, $y_{k-1} = q^k$ and $r_k = p_k/y_{k-1} = p$. For $\delta < 0$, we have $s_k = p_{k+\delta}/y_{k-1} = pq^{k+\delta}/q^k = pq^\delta$ when $k \geq -\delta$ and $s_k = 0$ otherwise. Also $\theta(k) = (k + \delta + 1)^+ pq^\delta$ and $m(n) = \lfloor -\log n/\log q \rfloor$. From Corollary 3.1, we obtain

$$(\log n)^{-1/2}(N_n + pq^\delta \log n/\log q) \xrightarrow[n]{D} N(0, -pq^\delta(q^{\delta+1} + q^\delta - 1)/\log q).$$

Weak records are observations such that $X_n \geq M_{n-1}$. In our context, they correspond to $\delta$-records with $\delta = -1$ and we have

$$(\log n)^{-1/2}(N_n + (p/q)\log n/\log q) \xrightarrow[n]{D} N(0, -(p/q^2)/\log q).$$

The above result was obtained by Bai *et al.* [3], using generating function methods. With some extra effort, our results could be extended to functional central limit theorems such as

$$(\log n)^{-1/2}(N_{\lfloor n^t \rfloor} + t(p/q)\log n/\log q) \xrightarrow[n]{D} \sqrt{-(p/q^2)/\log q}\, W(t)$$

for the number of weak records of geometric random variables. The limit $W(t)$ is the standard Wiener process and $\xrightarrow[n]{D}$ is understood as weak convergence on the Skorohod space $D[0,\infty)$.

For positive $\delta$, we apply Corollary 4.1 to obtain

$$(\log n)^{-1/2}(N_n + pq^\delta \log n/\log q) \xrightarrow[n]{D} N(0, -pq^\delta(q^{\delta+1} - (1 + 2\delta p)q^\delta + 1)/\log q).$$

**Example 5.2** *(Negative binomial)*. Here, $p_k = (-1)^k \binom{-a}{k} p^a q^k$ for $k \geq 0, 0 < p < 1, q = 1 - p$ and $a > 1$. From Vervaat ([22], Example 3.1), we have $p - (a-1)q/k \leq r_k \leq p$ and we obtain the same limiting distributions as the geometric example above.

**Example 5.3** *(Zeta)*. The zeta distribution has $p_k = (k+1)^{-a}/\zeta(a)$ for $k \in \mathbb{Z}_+$ and $a > 1$, with $\zeta(a) = \sum_{j=0}^{\infty}(j+1)^{-a}$. Here, $r_k = (k+1)^{-a}/\sum_{j=k}^{\infty}(j+1)^{-a} \underset{k}{\sim} (a-1)/k$. From Corollaries 3.1(b) and 4.1(b), we obtain

$$(\log n)^{-1/2}(N_n - \log n) \xrightarrow[n]{D} N(0,1).$$

Note that the normalizing sequences in this example do not depend on the value of $\delta$, positive or negative. This can be intuitively explained because samples from heavy-tailed distributions show, with high probability, values that are 'big' records.

**Example 5.4** *(Poisson)*. Let $p_k = e^{-\lambda}\lambda^k/k!$, $k \in \mathbb{Z}_+$, $\lambda > 0$. The following approximation of the failure rates $r_k$ can be found in Vervaat ([22], page 328):

$$\frac{\lambda}{k+1} - \left(\frac{\lambda}{k+1}\right)^2 \leq 1 - r_k \leq \frac{\lambda}{k+1}.$$

Let $\delta < 0$. Then it is easy to see that $\sum_{k=0}^{m(n)}(1-r_k)^{2\delta} \underset{n}{\sim} \lambda^{2\delta} \sum_{k=0}^{m(n)} k^{-2\delta} \underset{n}{\sim} \lambda^{2\delta} m(n)^{1-2\delta}/(1-2\delta)$ and $m(n)^{\delta-1/2}(\sum_{k=0}^{m(n)} s_k^\delta - \lambda^\delta m(n)^{1-\delta}/(1-\delta)) \xrightarrow[n]{} 0$, obtaining, from Theorem 3.1(b),

$$m(n)^{\delta-1/2}(N_n - \lambda^\delta(m(n))^{1-\delta}/(1-\delta)) \xrightarrow[n]{D} N(0, \lambda^{2\delta}/(1-2\delta)),$$

where $m(n) \underset{n}{\sim} \log n/\log\log n$.

When $\delta > 0$, we see, from Theorem 4.1(b), that the situation is quite different because given that $\sum_{k=1}^{n} e_k \underset{n}{\sim} \lambda^\delta \sum_{k=1}^{n} k^{-\delta}$, the number of $\delta$-records is finite if $\delta > 1$. For $\delta = 1$, we have $\sum_{k=1}^{n} e_k \underset{n}{\sim} \lambda \log n$ and it is easy to see that $\sum_{n=1}^{\infty} |\theta(n) - \lambda/(n+1)| < \infty$. Also, $m(n) \underset{n}{\sim} \log n/\log\log n$. Therefore,

$$(\log\log n)^{-1/2}(N_n - \lambda \log m(n)) \xrightarrow[n]{D} N(0, \lambda).$$

## 6. Concluding remarks

A referee suggested we consider the extension of our results to the case of $k$th upper order statistics, introducing the random quantity $S_{n,k} = \sum_{i=k+1}^{n} \mathbf{1}_{\{X_i > X_{i-1:i-k}+\delta\}}$, where $X_{i-1:i-k}$ denotes the $k$th upper order statistic of $X_1, \ldots, X_{i-1}$. It is easy to see that replacing $M_n = X_{n:n}$ by $X_{n:n-k+1}$ in (2.3) does not yield a martingale. However, the modification

$$S_{n,k} - \sum_{j=0}^{k-1} \theta(X_{n:n-j})$$

is a martingale. It is not clear, though, how to handle this process to get results analogous to those obtained in this paper.

# Appendix: Sums of minima and martingale central limit theorem

## A.1. Sums of partial minima

The martingale approach we use depends on asymptotic results for sums of partial minima of independent identically distributed random variables. The following weak law of large numbers from Deheuvels [7] is quite useful here.

Let $\{Z_n, n \geq 1\}$ be a sequence of independent identically distributed non-negative random variables, with common distribution function $G$, such that $G(z) > 0$ for all $z > 0$ and let $S_n = \sum_{i=1}^n \min\{Z_1, \ldots, Z_i\}$. Let also $G^-(t) = \inf\{z \geq 0 \mid G(z) \geq t\}$, for $0 \leq t < 1$ and $H(x) = \int_1^{e^x} G^-(1/u)\, du$ for $x \geq 0$.

**Theorem A.1** (Deheuvels [7], Theorem 7 and Corollary 4). *If $\lim_{x \to \infty} H(x)$ is finite, then $S_n$ grows almost surely to a finite limit as $n \to \infty$. Otherwise, if there is a sequence $x_n \uparrow \infty$ such that*

$$H(x_n + \log n)/H(\log n) \xrightarrow[n]{} 1 \tag{A.1}$$

*and*

$$\sum_{k=1}^n k G^-(1/k)^2 \Big/ \left(\sum_{k=1}^n G^-(1/k)\right)^2 \xrightarrow[n]{} 0, \tag{A.2}$$

*then*

$$S_n/H(\log n) \xrightarrow[n]{P} 1.$$

## A.2. A martingale central limit theorem

We use the martingale central limit theorem given by Hall and Heyde ([11], page 58), replacing the Lindeberg-type condition by the stronger Lyapunov-type condition (A.4).

**Theorem A.2.** *Let $\{\xi_i, i \geq 1\}$ such that $\mathrm{E}[|\xi_i|^3] < \infty$ and $\mathrm{E}[\xi_i | \mathcal{F}_{i-1}] = 0$, for all $i \geq 1$. For a sequence $b_n \uparrow \infty$, if the conditions*

$$\frac{1}{b_n^2} \sum_{i=1}^n \mathrm{E}[\xi_i^2 | \mathcal{F}_{i-1}] \xrightarrow[n]{P} 1 \tag{A.3}$$

and

$$\frac{1}{b_n^3} \sum_{i=1}^n \mathrm{E}[|\xi_i|^3|\mathcal{F}_{i-1}] \xrightarrow[n]{P} 0, \tag{A.4}$$

hold, then $\sum_{i=1}^n \xi_i/b_n \xrightarrow[n]{D} N(0,1)$.

**Lemma A.1.** *Let $\{a_n, n \geq 1\}$ be a sequence of positive terms such that $a_n \xrightarrow[n]{} \infty$ and $a_n/a_{n-1} \xrightarrow[n]{} 1$. Then $a_n/S_n \xrightarrow[n]{} 0$ and $S_{2,n}/(S_n)^2 \xrightarrow[n]{} 0$, where $S_n = \sum_{i=1}^n a_i$ and $S_{2,n} = \sum_{i=1}^n a_i^2$.*

**Proof.** The proof is a simple exercise. Let $\varepsilon > 0$ and take $N \in \mathbb{N}$ such that $a_n - a_{n-1} < \varepsilon a_n$ for all $n \geq N$. Then, for $n \geq N$,

$$a_n - a_0 = \sum_{i=1}^n (a_i - a_{i-1}) \leq a_N - a_0 + \varepsilon \sum_{i=N+1}^n a_i \leq a_N - a_0 + \varepsilon S_n,$$

which implies $a_n/S_n \xrightarrow[n]{} 0$. Analogously, if $a_n < \varepsilon S_n$ for $n > N$, then $S_{2,n} \leq S_{2,N} + \varepsilon \sum_{i=1}^n a_i S_i \leq S_{2,N} + \varepsilon (S_n)^2$, implying $S_{2,n}/(S_n)^2 \xrightarrow[n]{} 0$. □

**Lemma A.2** (Embrechts *et al.* [8], Proposition 3.1.1). *For $0 \leq \tau \leq \infty$ and a sequence $\{u_n, n \geq 1\}$, $n(1 - F(u_n)) \xrightarrow[n]{} \tau$ is equivalent to $P[M_n \leq u_n] \xrightarrow[n]{} e^{-\tau}$.*

**Corollary A.1.** *We have*

$$(\theta(M_n) - \theta(m(n)))/b_n \xrightarrow[n]{P} 0 \tag{A.5}$$

*if and only if $ny_{\Theta(\varepsilon b_n + \theta(m(n)))} \xrightarrow[n]{} 0$ and $ny_{\Theta(-\varepsilon b_n + \theta(m(n)))} \xrightarrow[n]{} \infty$ for all $\varepsilon > 0$.*

**Proof.** Convergence in (A.5) is equivalent to $P[\theta(M_n) \leq \varepsilon b_n + \theta(m(n))] \xrightarrow[n]{} 1$ and $P[\theta(M_n) \leq -\varepsilon b_n + \theta(m(n))] \xrightarrow[n]{} 0$ for all $\varepsilon > 0$. From Lemma A.2, these conditions are, respectively, equivalent to $nP[\theta(X_n) > \varepsilon b_n + \theta(m(n))] = ny_{\Theta(\varepsilon b_n + \theta(m(n)))} \xrightarrow[n]{} 0$ and $nP[\theta(X_n) > -\varepsilon b_n + \theta(m(n))] = ny_{\Theta(-\varepsilon b_n + \theta(m(n)))} \xrightarrow[n]{} \infty$. □

## Acknowledgements

The authors thank the referees for valuable comments and suggestions which greatly improved the presentation of the paper. Financial support is gratefully acknowledged from the FONDAP Project in Applied Mathematics, FONDECYT grants 1020836, 7020836 and 1060794, and MEC project MTM2004-01175. The authors are members of the research group Modelos Estocásticos (DGA).